\documentclass{amsart}

\usepackage{amsmath,amssymb,amsthm,latexsym}

\newtheorem{theorem}{Theorem}
\newcommand{\bt}{\begin{theorem}}
\newcommand{\et}{\end{theorem}}
\newtheorem{lemma}{Lemma}
\newcommand{\bl}{\begin{lemma}}
\newcommand{\el}{\end{lemma}}
\newtheorem{corollary}{Corollary}
\newcommand{\bc}{\begin{corollary}}
\newcommand{\ec}{\end{corollary}}
\newcommand{\bconj}{\begin{conjecture}}
\newcommand{\econj}{\end{conjecture}}
\newtheorem{problem}{Problem}
\newcommand{\bprob}{\begin{problem}}
\newcommand{\eprob}{\end{problem}}
\newcommand{\beq}{\begin{equation}}
\newcommand{\eeq}{\end{equation}}
\newcommand{\benum}{\begin{enumerate}}
\newcommand{\eenum}{\end{enumerate}}

\newcommand{\Z}{\ensuremath{\mathbf Z}}

\newcommand{\mca}{\ensuremath{ \mathcal A}}

\newcommand{\mcc}{\ensuremath{ \mathcal C}}
\newcommand{\mcd}{\ensuremath{ \mathcal D}}

\newcommand{\mbu}{\ensuremath{ \mathbf u}}

\newcommand{\bmat}{\left(\begin{matrix}}
\newcommand{\emat}{\end{matrix}\right)}
\newcommand{\bsmallmat}{\left(\begin{smallmatrix}}
\newcommand{\esmallmat}{\end{smallmatrix}\right)}

\DeclareMathOperator{\qqand}{\qquad\text{and}\qquad}

\title[Sizes of sumsets]{Inverse problems for sumset sizes of finite sets of integers}
\author{Melvyn B.  Nathanson}
\address{Department of Mathematics\\Lehman College (CUNY)\\Bronx, NY 10468}
\email{melvyn.nathanson@lehman.cuny.edu}

\date{\today}

\thanks{Supported in part by  PSC-CUNY Research Award Program grant 66197-00 54.}

\begin{document}

\begin{abstract}
Let $A$ be a finite set of integers and let $hA$ be its $h$-fold sumset.  
This paper investigates the sequence of sumset sizes  $( |hA| )_{h=1}^{\infty}$, 
the relations between these sequences for affinely inequivalent sets $A$ and $B$, 
and the comparative growth rates and configurations of the  sumset size 
sequences  $( |hA| )_{h=1}^{\infty}$ and $( |hB| )_{h=1}^{\infty}$. 
\end{abstract}

\maketitle

\section{The sequence of sumset sizes} 

Let $A$ be a nonempty  finite set of integers. 
The $h$-fold sumset of the set $A$, denoted  $hA$, is the set of all sums of $h$ 
not necessarily distinct elements of $A$.  
Let $|hA|$ denote the size (or cardinality) of the set  $|hA|$.  
A  \emph{direct problem} for sumset sizes in additive 
number theory is to compute the infinite sequence $( |hA| )_{h=1}^{\infty}$. 
For sufficiently large $h$, this sequence is completely determined 
by Theorem~\ref{sizes:theorem:Nathanson-2} and is eventually an arithmetic progression.  
For small $h$, Theorem~\ref{sizes:theorem:Lev} gives a lower bound for $|hA|$. 

Sets of integers $A$ and $B$ are  \emph{affinely equivalent} 
if there exist rational numbers $ \lambda$ and $\mu$  
with $\lambda \neq 0$ such that 
\[
B = \lambda\ast A + \mu = \{ \lambda x + \mu: x \in A\}.
\] 
For example, the sets $\{0,2\}$ and $\{1,4\}$ 
are affinely equivalent because 
\[
\{0,2\} = \frac{2}{3}\ast \{1,4\} - \frac{2}{3}.
\]
With $\lambda = -1$ and $a  \in \Z$, the set $A$ is affinely equivalent to 
its \emph{reflection}  $a  - A = (-1)\ast A + a  $.  
All 1-element sets of integers are affinely equivalent 
and all 2-element sets of integers are affinely equivalent.  
The sets $[0,k-1]$ and $[0,k-2] \cup \{k\}$ are affinely inequivalent  
sets of integers of size $k$ for all $k \geq 3$.  
If $A$ and $B$ are affinely equivalent sets with $B = \lambda\ast A + \mu$, then 
\[
hB = h(\lambda\ast A + \mu) = \lambda\ast hA + h\mu  
\]
and so $|hA| = |hB|$ for all $h \geq 1$.  
Thus, the size of a sumset is an \emph{affine invariant}.

It is natural to consider  \emph{inverse problems} 
for the sequence of sumset sizes $\left(|hA|\right)_{h=1}^{\infty}$ 
of  a nonempty  finite set $A$ of integers. 
To what extent does the sequence of positive integers 
$\left(|hA|\right)_{h=1}^{\infty}$ determine $A$?  
If $A$ and $B$ are affinely equivalent sets, then 
$\left(|hA|\right)_{h=1}^{\infty}=  \left(|hB|\right)_{h=1}^{\infty}$.  
What are the relations between the sequences 
$\left(|hA|\right)_{h=1}^{\infty}$ and $ \left(|hB|\right)_{h=1}^{\infty}$ 
for affinely inequivalent sets of the same size?  
Do there exist affinely inequivalent  finite sets of integers $A$ and $B$ 
such that $\left(|hA|\right)_{h=1}^{\infty}=  \left(|hB|\right)_{h=1}^{\infty}$? 
For every positive integer $h_1$, do there exists finite sets of integers $A$ and $B$ 
such that $|hA| = |hB|$ for all $h \leq h_1$ but  $|hA| < |hB|$ for all $h > h_1$? 
Do there exists finite sets of integers $A$ and $B$ 
such that $|hA| > |hB|$ for all $h \leq h_1$ but  $|hA| < |hB|$ for all $h > h_1$? 
We answer some of these question, but many  problems remain open.

{\emph{Notation.} In this paper, all sets are finite sets of integers.  
The \emph{interval of integers} $[u,v]$ is the set of all integers $n$ such that $u \leq n \leq v$. 
The \emph{integer part}, also called the \emph{floor}, of the real number $x$, denoted $\lfloor x \rfloor$, 
is the greatest integer $n \leq x$, and the  \emph{ceiling} of  the real number $x$, 
denoted $\lceil x \rceil$, is the smallest integer $n \geq x$.

\section{The direct problem for sumset sizes}
A finite set of integers is \emph{normalized} if $\min(A) = 0$ and $\gcd(A)=1$.
If the set $A$ is normalized, then the set $\max(A)-A$ is also normalized.  
Every finite set of at least  two integers  is affinely equivalent 
to an essentially unique  \emph{normalized} set. 
(If $A$ and $A'$ are affinely equivalent normalized sets of integers, then 
$A' = \max(A)-A$.) 

The following result is fundamental in additive number theory.

\bt[Nathanson~\cite{nath1972-7,nath96bb}]               \label{sizes:theorem:Nathanson}
Let $A$ be a finite set of $k \geq 2$ integers with $\min(A) = 0$, $\max(A) = a$, 
and   $\gcd(A)=1$. 
There exist integers $C$ and $D$ and subsets $\mcc \subseteq [0,C-2]$ 
and $\mcd \subseteq [0,D-2]$ such that 
\[
hA = \mcc \cup [C,ha-D]  \cup (ha- \mcd)
\]
for all $h \geq h_0(A)$. 
\et

The \emph{numerical semigroup} generated by a  set $A$ of relatively prime nonnegative integers 
is the set $S(A) = \{0\} \cup  \bigcup_{h=1}^{\infty} hA$.  This semigroup contains all sufficiently large integers.  
The \emph{Frobenius number} of the set $A$, denoted $g(A)$, 
is the largest integer that is not in the semigroup. 
The \emph{genus} of $A$, denoted $N(A)$, is  the number of nonnegative integers 
not contained in the numerical semigroup $S(A)$.   
(For the Frobenius problem, see Ram{\' i}rez-Alfons{\' i}n~\cite{rami05}.)
In Theorem~\ref{sizes:theorem:Nathanson}, 
\[
|\mcc| = C - N(A) \qqand |\mcd| = D - N(a-A)
\]
where $a-A = \{a-x:x\in A\}$ is the reflection of $A$.

Theorem~\ref{sizes:theorem:Nathanson} immediately implies the following.

\bt        \label{sizes:theorem:Nathanson-2}
Let $A$ be a  finite set of $k \geq 2$ integers 
with $\min(A) = 0$, $\max(A) = a$, and   $\gcd(A)=1$. 
Then 
\[
|hA| = ha + 1 - N(A) - N(a-A) 
\]
for all $h \geq h_0(A)$ and so the sequence of sumset sizes $\left(|hA|\right)_{h=1}^{\infty}$
is eventually an arithmetic progression with difference $a$.  
\et

The original result in~\cite{nath1972-7} was $ h_0(A) = a^2(k-1)$.  Recently, 
Andrew Granville, George Shakan, and Aled Walker~\cite{gran-shak20,gran-walk21} 
obtained $h_0(A) = a - k+2$.  
Vsevolod Lev~\cite{lev22} gives another proof of this result.  

The hard part of the direct problem is to understand $|hA|$ for small $h$.  
For this, the following result is useful.   

\bt[Lev~\cite{lev96b}]                           \label{sizes:theorem:Lev}
Let $A$ be a finite set of $k \geq 2$ integers with $\min(A) = 0$, $\max(A) = a$, 
and   $\gcd(A)=1$. 
For all $h \geq 2$, 
\[
|hA| \geq |(h-1)A| + \min(a, h(k-2)+1).
\]
\et

\section{Inverse problems for sumset size}

\bprob
Let $A$ be a normalized finite set of integers.  
The inverse problem for sumset sizes is to deduce  structural information 
about $A$ from the sumset size sequence $\left(|hA|\right)_{h=1}^{\infty}$.  
\eprob

For example, if $A$ is a normalized finite set of integers, 
then the sequence $\left(|hA|\right)_{h=1}^{\infty}$ is eventually 
an arithmetic progression with $|hA| = ha+1- N(A)-N(a-A)$ and so $\max(A) = a$. 
The set $A$ is an interval if and only if $2|hA| = |(h-1)A| + |(h+1)A|$
for all $h \geq 2$. 
The set $A$ is not an interval if $|hA| \neq h|A|-h+1$ for some $h \geq 2$.

The next result shows that the sumset size sequence 
$\left(|hA|\right)_{h=1}^{\infty}$  
does not determine the affine  equivalence class of $A$.

\bt      
For all $k \geq 3$, there exist affinely inequivalent sets $A$ and $B$ such that 
$|A| = |B| = k$ and $|hA| = |hB|$ for all $h \geq 2$.
\et

\begin{proof}
The proof is by explicit construction. 

Let $k=3$.  The sets 
\[
A = \{0,2,7\} \qqand B = \{0,3,7\}  
\]
are affinely inequivalent with $|A| = |B| = 3$.
We have 
\begin{align*}
2A & = \{0,2,4,7,9,14\} \\
2B & = \{0,3,6,7,10,14\} \\
3A & = \{ 0, 2, 4, 6, 7, 9, 11, 14, 16, 21 \} \\
3B & = \{ 0, 3, 6, 7, 9, 10, 13, 14, 17, 21 \} \\
4A & = \{0, 2, 4, 6, 7, 8, 9, 11, 13, 14, 16, 18, 21, 23, 28 \} \\
4B & = \{0, 3, 6, 7, 9, 10, 12, 13, 14, 16, 17, 20, 21, 24, 28\}
\end{align*}
and so 
\[
|2A| = |2B| = 6
\]
\[
|3A| = |3B| = 10 
\]
\[
|4A| = |4B| = 15.
\]
Applying Theorems~\ref{sizes:theorem:Nathanson} 
and~\ref{sizes:theorem:Nathanson-2},  we obtain finite sets  
\[
\mcc_A = \{0,2,4\} \qqand \mcd_A = \{ 0, 5, 7, 10, 12, 14, 15, 17, 19, 20, 21, 22 \}
\]
and 
\[
\mcc_B = \{0,3,6,7,9,10\} \qqand \mcd_B = \{0,4,7,8,11,12,14,15,16\}  
\]
such that 
\[
|\mcc_A | + |\mcd_A|  = |\mcc_B| + |\mcd_B| = 15.
\]
For all $h \geq 5$ we have 
\[
hA = \mcc_A \cup [6, 7h-24] \cup (7h - \mcd_A)
\]
and
\[
hB = \mcc_B \cup [12, 7h-18] \cup (7h - \mcd_B)
\]
and so  
\[
|hA| = |hB| = 7h -14.
\]
It follows that $|hA| = |hB|$ for all $h \geq 1$.

Let $k = 4$. 
The sets 
\[
A = \{0,3,5,6\} \qqand B = \{0,4,5,6\} 
\]
are affinely inequivalent with $|A| = |B| = 4$. 
For all $h \geq 2$ we have 
\[
hA = A \cup [8,6h]
\]
and
\[
hB = B\cup [8,6h]
\]
and so $|hA| = |hB| = 6h-3$ for all $h \geq 2$.

For all $k \geq 5$, the sets 
\[
A = \{0,3,5,6\}  \cup [8,k+3] \qqand B = \{0,4,5,6\}  \cup [8,k+3] 
\]
are affinely inequivalent with $|A| = |B| = k$. 
For all $h \geq 2$ we have 
\[
hA = A \cup [8,h(k+3)]
\]
and so $|hA| = |hB| = h(k+3)-3$ for all $h \geq 2$.
This completes the proof. 
\end{proof}

A useful experimental class of finite sets of integers consists of sets 
that are the union of a finite interval of integers and one additional element.  
The following result computes the sumset sizes of these sets.  
In Appendix~\ref{sizes:appendix},  there is a table of sumset sizes $|hA_{4,w}|$  for 
such sets of size 6, that is, sets of the form $A_{4,w} = \{0,1,2,3,4\} \cup \{w\}$ for $w \geq 5$.

\bt                                             \label{sizes:theorem:interval}
Let $h$,  $\ell$, and $w$ be positive integers with $w \geq {\ell}+1$ and let  
\[
A_{\ell,w} = [0,{\ell}] \cup \{w\}.  
\]
If $h \leq   (w-1)/\ell  $, then
\beq                                                     \label{sizes:interval-h-small}
|hA_{\ell,w}| = \frac{(h+1)({\ell}h+2)}{2}.
\eeq
If $h \geq  w/\ell $ and 
\[
 j_0 = j_0(h, \ell, w) = \left\lfloor h+1 - \frac{w}{{\ell}} \right\rfloor
 \]
then 
\beq                                                        \label{sizes:interval-h-big}
|hA_{\ell,w}| =    j_0w + 1  +  \frac{(h-j_0) ( 2 + {\ell}(h-j_0+1))}{2}.
\eeq 
\et

In Theorem~\ref{sizes:theorem:interval}, if $w = \ell q+r$ with $r \in [1,\ell]$, then 
$j_0 = h-q$.

\begin{proof}   
This is true for $h=1$ and all $\ell \geq 1$ because $|A_{\ell,w}| = \ell + 2$.  

Let $h \geq 2$.  For all $j \in [0,h]$, consider the interval of integers 
\[
I_j =  jw + (h-j)[0,{\ell}] = [jw, jw+ (h-j){\ell}].  
\]
Then 
\[
hA_{\ell,w}  = h([0,{\ell}] \cup \{w\}) = \bigcup_{j=0}^h \left( jw +  (h-j)[0,{\ell}]  \right) 
= \bigcup_{j=0}^h I_j .  
\]
For all $j \in [1,h]$, we have 
\[
\min(I_{j-1}) = (j-1)w < jw = \min(I_j) 
\]  
and (because $w \geq {\ell}+1$) 
\[
\max(I_{j-1}) = (j-1)w + (h-j+1){\ell}  < jw+ (h-j){\ell} = \max(I_j).
\]  
The  intervals $I_{j-1}$ and $I_j$ have a nonempty intersection if and only if 
\[
\min(I_j) \leq \max(I_{j-1}) 
\]
if and only if 
\[
jw \leq (j-1)w+(h-j+1){\ell} 
\]
 if and only if 
\[
w \leq  (h-j+1){\ell} 
\] 
if and only if 
\[
j \leq j_0 = \left\lfloor  h+1 - \frac{w}{{\ell}} \right\rfloor.
\]

If  $w \leq h{\ell} $, or, equivalently, if    $h \geq w/\ell$, then $j_0 \geq 1$, the $j_0+1$ intervals 
 $I_0, I_1,\ldots, I_{j_0}$ overlap, and 
\[
L =  \bigcup_{j=0}^{j_0} I_j = \left[ 0, j_0w+(h-j_0){\ell} \right]. 
\]
The $h+1-j_0$ intervals $L, I_{j_0+1}, \ldots, I_h$ are pairwise disjoint and so 
\begin{align*} 
|hA_{\ell,w}| =  \left| \bigcup_{j=0}^h I_j  \right| & =  \left|L \right|  
+  \sum_{j= j_0+1}^h  \left| I_j \right|  \\
& =   j_0w+(h-j_0){\ell}+1  + \sum_{j= j_0+1}^h ((h-j){\ell}+1) \\
& =   j_0w + 1  +  \frac{(h-j_0) ( 2 + {\ell}(h-j_0+1))}{2}.
\end{align*} 

If $w \geq h{\ell}+1$, or, equivalently, if $h \leq (w-1)/\ell$, 
then $j_0 = 0$ and  the intervals $I_0, I_1, \ldots, I_h$ are pairwise disjoint.  
We obtain  
\begin{align*} 
|hA_{\ell,w}| & = \sum_{j=0}^h \left| I_j \right|  = \sum_{j=0}^h \left(  (h-j){\ell}+1 \right) \\ 
& = \frac{(h+1)({\ell}h+2)}{2}.
\end{align*} 
This completes the proof.  
\end{proof}

\bt               \label{sizes:theorem:A1w}
Let $h$  and $w$ be positive integers with $w \geq 2$ and let  
\[
A_{1,w} = \{0, 1,w\}.  
\]
If  $h \leq w-1$,   then 
\beq                                                     \label{sizes:A1w-small}
|hA_{1,w}| = \frac{(h+1)(h+2)}{2}.
\eeq
If $h \geq w $, then 
\beq                                                        \label{sizes:A1w-big}
|hA_{1,w}|  = hw + 1 - \frac{(w-1)(w-2)}{2}.
\eeq 
\et

\begin{proof}
We apply Theorem~\ref{sizes:theorem:interval} 
with $\ell = 1$ and $j_0 = h+1-w$.  
If $h \leq w-1$, then we immediately obtain~\eqref{sizes:A1w-small}.
If $h \geq w$, then $h-j_0 = w-1$ and 
\begin{align*}
|hA_{1,w}| & = (h+1-w)w + 1 + \frac{(w-1)(w+2) }{2} \\
& = hw + 1 - \frac{(w-1)(w-2)}{2}.
\end{align*}
This completes the proof.  
\end{proof}

To obtain formula~\eqref{sizes:A1w-big}, 
we can also use  the genus formula of Sylvester~\cite{sylv82a}: 
If $v$ and $w$ are relatively prime positive integers and $A = \{0,v,w\}$, then 
\[
N(A) = \frac{(v-1)(w-1)}{2}. 
\] 
Thus,  
\[
N(A_{1,w} ) = N(\{0,1,w\}) = 0
\]
 and 
\[
N(w-A_{1,w} ) = N(\{0,w-1,w\}) =  \frac{(w-1)(w-2)}{2}.
\]
  Applying Theorem~\ref{sizes:theorem:Nathanson-2} 
with $h \geq h_0(A_{1,w} ) = w-1$, we obtain 
\[
|hA_{1,w} | = hw+1-\frac{(w-1)(w-2)}{2}. 
\]

\bt                                             \label{sizes:theorem:interval-A}
Let $h$,  $\ell$, and $w$ be positive integers with $w \geq {\ell}+1$.  
Let $A'$ be a nonempty subset of $[0,\ell]$ and let 
\[
A  = A' \cup \{w\}.  
\]
If $h \geq w/\ell$ and 
\[
 j_0 = \left[  h+1 - \frac{w}{{\ell}} \right] 
 \]
then 
\[
|hA | \leq  j_0w + 1  +  \frac{(h-j_0) ( 2 + {\ell}(h-j_0+1))}{2}.
\] 
\et

\begin{proof}
Let 
\[
A_{\ell,w} = [0,{\ell}] \cup \{w\}.  
\]
We have $A \subseteq A_{\ell,w}$.   
If $h \geq w/\ell$, then 
\[
|hA| \leq |hA_{\ell,w}| =    j_0w + 1  +  \frac{(h-j_0) ( 2 + {\ell}(h-j_0+1))}{2}.
\]
This completes the proof.  
\end{proof}

\section{Sumset size oscillations}

Consider the following ``experimental'' data from Appendix~\ref{sizes:appendix} 
for sumsets of sets of the form $A_{4,w} = [0,4] \cup \{w\}$.

We have $|2A_{4,9}| = |2A_{4,10}| = 15$ and 
$|hA_{4,10}| - |hA_{4,9}| = h-2 > 0$ for all $h \geq 3$.
 
We have $|2A_{4,13}| = |2A_{4,14}| = 15$, 
$|3A_{4,13}| = |3A_{4,14}| = 28$, and 
$|hA_{4,14}| - |hA_{4,13}| = h-3> 0$ for all $h \geq 4$.
 
We have $|2A_{4,17}| = |2A_{4,18}| = 15$, 
$|3A_{4,17}| = |3A_{4,18}| = 28$, $|4A_{4,17}| = |4A_{4,18}| = 45$, and 
$|hA_{4,18}| - |hA_{4,17}| = h-4> 0$ for all $h \geq 5$.

We have $|2A_{4,21}| = |2A_{4,22}| = 15$, 
$|3A_{4,21}| = |3A_{4,22}| = 28$, $|4A_{4,21}| = |4A_{4,22}| = 45$, 
$|5A_{4,21}| = |5A_{4,22}| = 66$, and 
$|hA_{4,22}| - |hA_{4,21}| = h-5> 0$ for all $h \geq 6$.

These examples suggest the following result.  

\bt                               \label{sizes:theorem:BigSplit}
For all  integers $k \geq 3$ and $h_1 \geq 1$, there exist finite sets $A$ and $B$ 
with $|A| = |B| = k$ such that $|hB| = |hA|$ for all $h \leq h_1$ and 
$ |hB| = |hA| + h-h_1 > |hA|$ for all $h \geq h_1 +1$.
\et

\begin{proof}
For integers $k \geq 3$ and $h_1 \geq 1$, define  
\[
\ell = k-2 \qqand w =  \ell (h_1+1). 
\] 
We apply Theorem~\ref{sizes:theorem:interval} to the sets 
\[
A = A_{\ell,w} = [0,\ell] \cup \{w\}
\]
and 
\[
B = A_{\ell,w+1} = [0,\ell] \cup \{w+1\}. 
\]

First, we compute the sumset size $|hA|$. 
For all 
\[
h \leq \left\lfloor \frac{w-1}{\ell} \right\rfloor =  \left\lfloor \frac{ \ell (h_1+1) -1}{\ell} \right\rfloor =  \left\lfloor h_1 +1-\frac{1}{\ell} \right\rfloor = h_1
\]
formula~\eqref{sizes:interval-h-small} of Theorem~\ref{sizes:theorem:interval} gives 
\[
|hA| = \frac{(h +1)(\ell h+2)}{2}. 
\]
For all
\[
h \geq \left\lceil \frac{w}{\ell} \right\rceil =  \left\lceil \frac{\ell(h_1+1)}{\ell} \right\rceil =  h_1+1 
\]
we have 
\[
j_0(h,\ell,w) =  \left\lfloor h+1-\frac{w}{\ell}  \right\rfloor 
=  \left\lfloor  h+1-\frac{ \ell (h_1+1) }{\ell}  \right\rfloor  = h-h_1. 
\]
From formula~\eqref{sizes:interval-h-big} of Theorem~\ref{sizes:theorem:interval}, 
we obtain
\begin{align*}
|hA| & = (h-h_1)w+1+\frac{h_1(2+\ell(h_1+1))}{2} \\
& = (h-h_1)w+1+\frac{h_1(w+2)}{2}.
\end{align*}
In particular, for $h = h_1+1$, we obtain 
\begin{align*}
|(h_1+1)A| &  =  w + 1  +  \frac{h_1  (w+ 2) }{2}. 
\end{align*}

Next, we compute the sumset size $|hB|$. 
For all 
\[
h \leq \left\lfloor \frac{(w+1)-1}{\ell} \right\rfloor =  \left\lfloor \frac{w}{\ell} \right\rfloor 
= \left\lfloor \frac{\ell (h_1+1)}{\ell} \right\rfloor = h_1+1
\]
formula~\eqref{sizes:interval-h-small} of Theorem~\ref{sizes:theorem:interval} gives 
\[
 |hB| = \frac{(h +1)(\ell h+2)}{2} 
\]
and so $|hB| = |hA|$ for all $h \leq h_1$. 
In particular, for $h = h_1+1$, we obtain 
\begin{align*} 
|(h_1+1)B|  & = \frac{(h_1+2)({\ell}(h_1+1)+2)}{2}  \\
& = \frac{(h_1+2) (w+2)}{2}  \\ 
 & =  w+2+\frac{h_1(w+2)}{2} \\ 
 & =   |(h_1 +1)A| +1.  
\end{align*}
For all 
\[
h \geq \left\lceil \frac{w+1}{\ell} \right\rceil =  \left\lceil \frac{\ell(h_1+1)+1}{\ell} \right\rceil =  h_1+2
\]
we have 
\begin{align*}
j_0 = j_0(h,\ell,w+1) &  = \left\lfloor h+1-\frac{w+1}{\ell} \right\rfloor \\
& = \left\lfloor h+1-\frac{ \ell  (h_1+1) +1}{\ell}\right\rfloor \\
& = \left\lfloor h - h_1-\frac{1}{\ell} \right\rfloor = h-h_1-1.
\end{align*}
From formula~\eqref{sizes:interval-h-big}, we obtain
\begin{align*}
|hB| & = (h-h_1-1)(w+1)+1+\frac{(h_1+1)(2+\ell(h_1+2))}{2} \\
& = (h-h_1)w   + (h-h_1) - w + (h_1+1) + \frac{w(h_1+2)}{2}\\ 
& =   (h-h_1)w   + 1  + \frac{h_1(w+2)}{2}  + h-h_1 \\
& = |hA| + h-h_1.
\end{align*}
This completes the proof. 
\end{proof}

\section{Open problems}
The comparative study of sumset sizes of finite sets of integers is a new research 
topic in additive number theory with many natural open problems.   

We begin with an example of a pair of sets with oscillating sumset sizes.  
Let  $ h_2, g,  \ell, b$, and  $w$ be positive integers such that  
\[
2 \leq h_2 < g  < g^{\ell} < b < w  \qqand h_2\ell < w.
\] 
Consider the sets 
\[
A _0 = [0,\ell]    \qqand 
G_0 = \left\{  1,g,g^2,\ldots, g^{\ell -1} \right\}  
\]
and 
\[
A = A_{\ell,w} =A_0  \cup \{ w \}   \qqand 
G = G_0  \cup \{ b\}.
\]
For $h_1 = 1$, we have 
\beq                   \label{sizes:1}
|h_1A| > |h_1B|.
\eeq    

For all $j \in [0,h_2]$, we have 
\[
I_j = jA_0 +(h_2-j)w = [(h_2-j)w, (h_2-j)w + j\ell]. 
\] 
For  all $j \in [1,h_2]$, the inequality $h_2\ell < w$ implies 
\[
(h_2-j)w + j\ell < (h_2-(j-1))w 
\]
and so the $h_2+1$  intervals $I_0, I_1, \ldots, I_{h_2}$ are pairwise disjoint.  
From  
\[ 
h_2A = \bigcup_{j=0}^{h_2} \left( jA_0 +(h_2-j)w \right) =  \bigcup_{j=0}^{h_2} I_j
\]
we obtain 
\[
\left| h_2A \right| = \sum_{j=0}^{h_2} |I_j| = \sum_{j=0}^{h_2}  |jA_0|. 
\]

Because $h_2 < g$ and $h_2g^{\ell-1} < b$, the set $G_0$ is a $B_{h_2}$-set, that is, a set such that 
all sums of at most $h_2$ not necessarily distinct elements are distinct, 
and  the set $G = G_0 \cup \{b\}$ is also  a $B_{h_2}$-set. 
It follows that 
\[
\left| jG_0 + (h_2-j) b \right| = |jG_0| = \binom{ \ell + j -1}{j} 
\]
and the $h_2+1$ sets $jG_0 +(h_2-j)  $ are also pairwise disjoint for $j \in [0,h_2]$. 
From 
\[
h_2B = \bigcup_{j=0}^{h_2} \left( jG_0 + (h_2-j) b \right)   
\]
we obtain
\[
\left| h_2 B \right| = \sum_{j=0}^{h_2} \left| jG_0 + (h_2-j) b \right| 
= \sum_{j=0}^{h_2} \left| jG_0 \right|. 
\] 
For $j \in [2,h_2]$, we have 
\[
 \left| jA_0 \right| = j\ell + 1 <  \binom{ \ell + j -1}{j} =  \left|  jG_0 \right|  
\]
and so  
\beq                   \label{sizes:2}
|h_2A| < |h_2B|. 
\eeq

By Theorem~\ref{sizes:theorem:Nathanson-2}, the choice of 
$b < w$ implies that there is an integer $h_3 > h_2$ such that 
$ |hA| > |hB|$ for all $h \geq h_3$.  In particular, we have positive integers 
\beq                   \label{sizes:3}
|h_3A| > |h_3B|. 
\eeq
Comparison of inequalities~\eqref{sizes:1},~\eqref{sizes:2}, and~\eqref{sizes:3}  
 suggests the following questions.

\bprob 
For every  integer $m \geq 3$, do there exist finite sets $A$ and $B$ 
and an increasing sequence of positive integers $h_1 < h_2 < \cdots < h_m$ such that 
\[
|h_i A|  > |h_i  B|   \qquad \text{if $i$ is odd} 
\]
and 
\[
  |h_i A| <  |h_i  B|   \qquad \text{if $i$ is even.} 
\]
Do there exist such sets with $|A| = |B|$? 
Do there exist such sets with $|A| = |B|$ and $\max(A) = \max(B)$?
\eprob

\bprob
Let  $(h_i)_{i=0}^m$ be a increasing sequence of positive integers.  
Do there exist  finite sets $A$ and $B$ such that 
\[
|h_iA| >  |h_i  B|   \qquad \text{if $i$ is odd} 
\]
and 
\[
 |h_iA| <  |h_i  B|    \qquad \text{if $i$ is even.} 
\]
Do there exist such sets with $|A| = |B|$?   
Do there exist such sets with $|A| = |B|$ and $\max(A) = \max(B)$? 
\eprob

\bprob
Prove or disprove the following statement:  For every positive integer $h^*$ 
there exist finite sets of integers $A$ and $B$ such that 
\[
|hA|  > |hB| \qquad \text{if $h$ is odd and $h \leq h^*$} \\
\]
and 
\[
|hA| <  |hB|  \qquad \text{if $h$ is even and $h \leq h^*$.}
\]
\eprob

\bprob
 Consider the analogous oscillation problems for sums of finite sets of lattice points in $\Z^n$
 and also for restricted sumsets, that is, sumsets of the form 
\[
\widehat{hA} = \left\{ a_1+\cdots  +a_h: a_i \in A \text{ and } a_i \neq a_j \text{ for } i \neq j\right\}.
\]
\eprob

\section{$\tau$-type problems}
Next we consider a more general class of oscillation  problems for $n$-tuples of 
sumset sizes of finite sets.  

The \emph{normalization} of an $n$-tuple $\mbu = (u_1,\ldots, u_n)$ 
of integers is the $n$-tuple $\tau(\mbu) = (\tau(1), \tau(2), \ldots, \tau(n))$ obtained 
by replacing the $i$th smallest entry in $\mbu$ with $i$ for all $i$.  
A \emph{$\tau$-tuple of length $n$} is an $n$-tuple of the form $\tau(\mbu)$ 
for some $n$-tuple \mbu.

Here are examples of $6$-tuples $\mbu$ and their associated $\tau(\mbu)$:
\vspace{0.3cm} 
\begin{center} 
\begin{tabular}{  l | c | c| c| c }  
$\mbu$         &    (1,2,3,4,5,6)   &  (-2,13,11,0,22,4)  &  (7,3,2,9,3,5)  & (9,7,8,9,7,8)  \\  \hline
$\tau(\mbu)$ &   (1,2,3,4,5,6)   & (1,5,4,2,6,3)        & (4,2,1,5,2,3) & (3,1,2,3,1,2) 
\end{tabular}
\end{center}
\vspace{0.3cm} 
We have $\tau(\sigma(1),\sigma(2),\ldots, \sigma(n) ) = (\sigma(1),\sigma(2),\ldots, \sigma(n) )$ 
for all permutations $\sigma$ of $[1,n]$.  

Let $\mca = (A_1, A_2, \ldots, A_n)$ be an $n$-tuple of finite sets. 
Associated to the $n$-tuple of sets $h$-fold sumsets $h\mca = (hA_1, hA_2, \ldots, hA_n)$ 
is the $n$-tuple of integers $|h\mca| =  ( |hA_1|, |hA_2|, \ldots, |hA_n|)$ 
and its normalization $\tau(|h\mca|)$.

\bprob
Let $\mbu_1, \mbu_2, \ldots, \mbu_m$ be a finite sequence of $\tau$-tuples 
of integers  of length $n$.
Does  there exist an $n$-tuple of  finite sets $\mca = (A_1, A_2, \ldots, A_n)$  
and a finite increasing sequence of positive integers 
$h_1 < h_2 < \cdots < h_m$ such that 
\[
\tau(|h_k \mca|) = \tau(|h_kA_1|, |h_kA_2|, \ldots, |h_kA_n|) = \tau(\mbu_k)
\]
for all $k \in [1,m]$?
\eprob

For example, let 
\[
\mbu_1 = (1,2,3), \qquad \mbu_2 = (2,3,1), \qquad \mbu_3 = (3,1,2).
\] 
Does  there exist a triple of  finite sets $(A_1, A_2, A_3)$ and integers $h_1 < h_2 < h_3$ 
such that 
\[
\tau(|h_1 \mca|) = \tau(|h_1A_1|, |h_1A_2|,  |h_1A_3|) =  (1,2,3)
\]
\[
\tau(|h_2 \mca|) = \tau(|h_2A_1|, |h_2A_2|,  |h_2A_3|) =  (2,3,1)
\]
\[
\tau(|h_3 \mca|) = \tau(|h_3A_1|, |h_3A_2|,  |h_3A_3|) =  (3,1,2)
\]
or, equivalently, 
\[
|h_1A_1| < |h_1A_2| < |h_1A_3|
\]
\[
|h_2A_3|  < |h_2A_1| < |h_2A_2|  
\]
\[
|h_3A_2| < |h_3A_3| < |h_3A_1|. 
\]

\bprob
Let $\mbu_0, \mbu_1, \mbu_2, \ldots, \mbu_m$ be a finite sequence 
of $n$-tuples of integers.
Does  there exist an $n$-tuple of  finite sets $(A_1, A_2, \ldots, A_n)$  
and a finite increasing sequence of positive integers 
$h_1 < h_2 < \cdots < h_m$ such that 
\[
\tau(|\mca |) = \mbu_0
\]
and 
\[
\tau(|h_k \mca|) = \tau(|h_kA_1|, |h_kA_2|, \ldots, |h_kA_n|) = \tau(\mbu_k)
\]
for all $k \in [1,m]$?
\eprob

\bprob
Let $\mbu_1, \mbu_2, \ldots, \mbu_m$ be a finite sequence 
of $n$-tuples of integers.  
Can one compute an integer $N = N(\mbu_1, \mbu_2, \ldots, \mbu_m)$ with the following property:
There exist an $n$-tuple of  finite sets $\mca = (A_1, A_2, \ldots, A_n)$  
and a finite increasing sequence of positive integers 
$h_1 < h_2 < \cdots < h_m$ such that 
\[
\tau(|h_k \mca|)  = \tau(|h_kA_1|, |h_kA_2|, \ldots, |h_kA_n|) = \tau(\mbu_k)
\]
for all $k \in [1,m]$ 
if and only if there exist such sets with $A_i \subseteq [0,N]$ for all $i \in [1,n]$?
\eprob

\section{Added January 31, 2025}
Noah Kravitz~\cite{krav25} 
recently proved the following theorem, which answers some of the questions 
posed in this paper. 

\bt
Let $n$ and $m$ be positive integers and let $\mbu_1, \mbu_2, \ldots, \mbu_m, \mbu_{\infty}$ 
be a finite sequence of $n$-tuples of integers.  
There is an $n$-tuple of  finite sets $\mca = (A_1, A_2, \ldots, A_n)$  such that 
\[
\tau(|h  \mca|)  = \tau(|h A_1|, |h A_2|, \ldots, |h A_n|) = \tau(\mbu_h)
\]
for all $h \in [1,m]$ and 
\[
\tau(|h  \mca|)  = \tau(|h A_1|, |h_kA_2|, \ldots, |h_kA_n|) = \tau(\mbu_{\infty})
\]
for all $h > m$. 
\et

\newpage

\appendix    
\section{Table of $|hA_{4,w}|$ for $A_{4,w} = [0,4] \cup \{w \}$, \\ 
$h \in [2,9]$ 
and $w \in [5,25] \cup \{30,35, 40,45, 50\}$}   \label{sizes:appendix} 

\ \\
\begin{center} 
\begin{tabular}{ | r | r  | r  | r  | r  | r  | r  | r  | r  |}  
\hline
$|hA_{4,w}|$ & h=2 & 3 & 4 & 5 & 6 & 7 & 8 & 9 \\ \hline
w=5  & 11 & 16 & 21 & 26 & 31 & 36 & 41 & 46  \\ \hline 
6 &  12 & 18 & 24 & 30 & 36 &  42 & 48 & 54  \\ \hline 
7 & 13 & 20 & 27 & 34 & 41 & 48 & 55 & 62  \\ \hline 
8 & 14 & 22 & 30 & 38 & 46 & 54 & 62 & 70 \\ \hline 
9 & 15 & 24 & 33 & 42 & 51 & 60 & 69 & 78  \\ \hline 
10 & 15 & 25 & 35 &45 & 55 & 65 & 75 & 85  \\ \hline 
11 & 15 & 26 & 37 & 48 & 59 &70 & 81 & 92 \\ \hline
12 & 15 & 27 & 39 &  51 & 63 & 75 & 87& 99 \\  \hline  
13 & 15 & 28 & 41 &   54 & 67 & 80 & 93 & 106  \\ \hline 
14 & 15 & 28 & 42 & 56 & 70 & 84 & 98 & 112 \\ \hline 
15 & 15 & 28 & 43 & 58 & 73 & 88 & 103 & 118  \\ \hline 
 16  & 15 & 28 & 44 & 60 & 76 & 92  & 108  & 124 \\ \hline 
  17 &15 & 28 & 45 & 62 & 79 & 96 & 113 &  130 \\ \hline 
18   & 15 & 28 & 45 & 63 & 81 & 99 & 117 & 135 \\ \hline 
  19  & 15 & 28 & 45 & 64 & 83 & 102 & 121 & 140  \\ \hline 
20   & 15 & 28 & 45 & 65 & 85 & 105 & 125  & 145 \\ \hline 
  21  &15 & 28 & 45 & 66 & 87 & 108 & 129 & 150 \\ \hline 
22  &15 & 28 & 45 & 66 & 88 & 110 & 132 & 154 \\ \hline 
23  &15 & 28 & 45 & 66 & 89 & 112 & 135 & 158 \\ \hline 
24  &15 & 28 & 45 & 66 & 90 & 114  & 138 & 162 \\ \hline 
25  &15 & 28 & 45 & 66 & 91 & 116  & 141 & 166  \\ \hline 
30  & 15 & 28 & 45 & 66 & 91 & 120 & 150  & 180 \\ \hline 
35 & 15 & 28 & 45 & 66 & 91 & 120 & 153 &188  \\ \hline
40  & 15 & 28 & 45 & 66 & 91 & 120 & 153  & 190 \\ \hline 
45 & 15 & 28 & 45 & 66 & 91 & 120 & 153 & 190 \\ \hline
50  & 15 & 28 & 45 & 66 & 91 & 120 & 153  & 190 \\ \hline 
\end{tabular}
\end{center}

\end{document}